\documentclass[12pt]{article}
\usepackage{graphicx}
\usepackage{latexsym,amssymb,amsmath,amscd,amsfonts,mathabx}
\topmargin - 2cm
\oddsidemargin - 0.1cm
\evensidemargin - 5cm
\newtheorem{theorem}{Theorem}[section]

\newtheorem{lemma}[theorem]{Lemma}

\newtheorem{definition}[theorem]{Definition}
\newtheorem{remark}[theorem]{Remark}

\newcommand{\R}{\mathbb{R}}
\newcommand{\C}{\mathbb{C}}

\begin{document}

\title{\bf On fractional parabolic systems of vector order}

\author{Ravshan Ashurov$^{1,2}$,  \  Ilyoskhuja Sulaymonov$^3$}
\date{}
\maketitle

\begin{center}
{\it $^1$V.I. Romanovskiy Institute of Mathematics, Uzbekistan Academy of Science, Tashkent, 100174, Uzbekistan;\\

$^2$ University of Tashkent for Applied Sciences, Str. Gavhar 1, Tashkent 100149, Uzbekistan; e-mail: ashurovr@gmail.com\\

$^3$ National University of Uzbekistan, Tashkent, Uzbekistan; e-mail: ilyosxojasulaymonov@gmail.com}

\end{center}
\vspace{10pt}

\begin{abstract}
The paper considers the Cauchy problem for the system of partial differential equations of fractional order $D_t^{\mathcal{B}} {U}(t,x) + \mathbb{A}(D) {U} (t,x)=H(t,x) $. Here $U$ and $H$ are vector-functions, the $m\times m$ matrix of differential operators $\mathbb{A}(D)$ is triangular (elements above or below the diagonal are zero). Operators located on the diagonal are elliptic. The main distinctive feature of this system is that the vector-order $\mathcal{B}$ has different components $\beta_j\in (0,1]$, and $\beta_j$ are not necessarily rational.  Sufficient conditions (in some cases they are necessary) on the initial function and the right-hand side of the equation are found to ensure the existence of a classical solution. Note that the existence of a classical solution to systems of fractional differential equations was studied by various authors, but in all these works the fractional order had the same components for each equation: $\beta_j=\beta$, $j=1,...,m$.
\end{abstract}

%
\vspace{2pc}

%
%
%
%


\section{Introduction}
In recent years, there has been great interest among specialists in the study of systems of fractional differential equations (ordinary and partial differential equations), primarily in connection with the rich applications of such systems in modeling various processes (see applications in biosystems \cite{DasGupta}, \cite{Rihan}, in ecology \cite{Khan}, \cite{Rana}, in epidemiology \cite{Zeb}, \cite{Islam}, \cite{Lyapunov-2},  \cite{Almeida}, \cite{Rajagopal}, etc).

When studying systems of differential equations, as well as in their applications, an important place is occupied by establishing a representation of the solution of the system. That is why the research of many authors is devoted to the search for such a formula for a solution of various systems of equations.
For systems of ordinary differential equations of fractional order, finding a formula for the solution is relatively simple (see, for example, \cite{Varsha} and  \cite{Veber}).

Recently, S. Umarov published a number of fundamental works \cite{UACh, Umarovappear, UmarovNew}, where the representation of solutions for various systems of fractional differential equations was established. The most general systems of differential equations are considered in work \cite{UmarovNew}, the results of which are new even for systems of ordinary fractional differential equations. The author, in particular, established a formula for a solution of the following Cauchy problem
\begin{align}\label{Cauchy_01_U}
  	D_t^{\mathcal{B}} {U}(t,x) &+ \mathbb{A}(D) {U} (t,x)=H(t,x), \quad t>0, \ x \in \mathbb{R}^n,
  	\\
  	U(0,x) & =\varPhi(x), \quad x \in \mathbb{R}^n, \label{Cauchy_02_U}
  \end{align}
where $(t, x)\in \mathbb{R}_+\times \mathbb{R}^n$, $\mathcal{B}=\langle \beta_1,\beta_2,\dots,\beta_m\rangle$, $\beta_j\in(0,1], j=1\dots m$, $\mathbb{A}(D)=\{A_{i,j}(D)\}$ is a $m\times m$ matrix whose elements are pseudo-differential operators,
$$
\Phi(x)=\langle \varphi_1(x), \varphi_2(x), \dots, \varphi_m(x)\rangle,
$$
$$
H(t,x)=\langle h_1(t,x), h_2(t,x), \dots, h_m(t,x)\rangle,
$$
are a given vector-functions, 
$$
U(t,x)=\langle u_1(t,x), u_2(t,x), \dots, u_m(t,x)\rangle,
$$
is an unknown vector-function and 
$$
D_t^{\mathcal{B}}U(t,x)=\langle D_t^{\beta_1}u_1(t,x), D_t^{\beta_2}u_2(t,x), \dots, D_t^{\beta_m}u_m(t,x)\rangle,
$$
is the fractional order derivative of order $ 0< \beta_j \leq 1, j=1, \dots, m$ in the sense of Caputo, which for any continuous function $f: \mathbb{R}_+\to \mathbb{R}$ is defined as (see, for example \cite{KLB}, p. 91): 
\[
D^{\beta_j}_t f(t)= \frac{1}{\Gamma(1-\beta_j)}\frac{d}{dt} \int\limits_0^t \frac{f(\xi)-f(0)}{(t-\xi)^{\beta_j}}d\xi, \,\,\, t>0,
\]
provided that the right side is pointwise defined on $\mathbb{R}_+$.

What is important here is the fact that the components $\beta_j$ of the vector order $\mathcal{B}$ are arbitrary real numbers from the interval $(0,1]$. The fact is that if all $\beta_j$ are the same: $\beta_j=\beta, \, j=1,..., m$, then the representation for the solution is quite simple (see the review work A.N. Kochubei \cite{Kochuei1} and \cite{Umarovappear}), and if $\beta_j$ are rational numbers, then this case can be reduced to the previous case in the standard way (see \cite{Umarovappear}). 

Let us note an important feature of the work of S. Umarov \cite{UmarovNew}. The solution $U(t,x)$ to problem (\ref{Cauchy_01_U}), (\ref{Cauchy_02_U}) is found in the class $\Psi$ of functions whose Fourier transform has a compact support. The initial function $\Phi(x)$ and the right-hand side of the equation $H(t,x)$ are also taken from the same class $\Psi$. The advantage of this choice is that when applying the Fourier transform to equation (\ref{Cauchy_01_U}) (and the method of work \cite{UmarovNew} is based precisely on the sequential application of the Fourier and Laplace transforms), their existence is beyond doubt. But on the other hand, the class of functions $\Psi$ is narrow, and since the class of pseudo-differential matrix operators $\mathbb{A}(D)$  in work \cite{UmarovNew} is quite wide, it is impossible to obtain the corresponding estimates for the solution, which would make it possible to close the class $\Psi$.

The question naturally arises whether it is possible to construct a solution to the Cauchy problem (\ref{Cauchy_01_U}), (\ref{Cauchy_02_U}) in ordinary classical classes, for example, a classical solution. It is to the solution of this issue that the work \cite{AO} is devoted, where, considering the order of fractional derivatives $\beta_j=\beta$, $j=1, ..., m$, the authors obtained sufficient (in some cases they are also necessary) conditions for the initial function $\Phi$ and for the right-hand side of equation $H(t,x)$, which guarantee the existence of a classical solution to the Cauchy problem (\ref{Cauchy_01_U}), (\ref{Cauchy_02_U}).

The purpose of this work is to establish a similar result in the case when the vector of the order of fractional derivatives $\mathcal{B}$  has arbitrary components $\beta_j\in (0,1]$. To do this, consider a lower triangular matrix of an arbitrary differential expressions $A_{i,j}(D)=\sum\limits_{|\alpha|\leq \ell_{i,j}}a_{\alpha} D^\alpha$ with the order $\ell_{i,j}$ and constant coefficients:
\[
\mathbb{A}(D)= \begin{bmatrix}
	    A_{1,1}(D)& 0 & \dots & 0 \\
	    A_{2,1}(D)& A_{2,2}(D) &\dots & 0\\
	    \dots  & \dots  & \dots & \dots  \\
	    A_{m, 1}(D) & A_{m,2}(D) & \dots & A_{m,m} (D)
	    \end{bmatrix},\] 
with the matrix-symbol
\[
\mathcal{A}(\xi) = \begin{bmatrix}
	    A_{1,1}(\xi)& 0 & \dots & 0 \\
	    A_{2,1}(\xi)& A_{2,2}(\xi) &\dots & 0\\
	    \dots  & \dots  & \dots & \dots  \\
	    A_{m, 1}(\xi) & A_{m,2}(\xi) & \dots & A_{m,m} (\xi)
	    \end{bmatrix},
\]
where  $\alpha=(\alpha_1, \alpha_2, \dots, \alpha_n)$ - multi-index and $D=(D_1, D_2, \dots, D_n)$,
$D_j=\frac{1}{i}\frac{\partial}{\partial x_j}$, $i=\sqrt{-1}$ and $A_{i,j}(\xi)=\sum\limits_{|\alpha|\leq \ell_{i,j}}a_{\alpha} \xi^\alpha$ is a symbol of operator $A_{i,j}(D)$.

We will further assume that all diagonal operators $A_{j,j}(D)$, $j=1, \dots, m$, are homogeneous elliptic differential operators and have the highest order among the operators acting on the function $u_j(t, x)$: 
\begin{equation}\label{maxorder}
    \ell_{j,j}>\ell_{i, j}, \quad \text{for all} \quad i\neq j, \quad i,j=1,\dots, m.
\end{equation}
By $p^\ast$ and $q_j$ we denote the following expressions:
\begin{equation}\label{p}
 p^\ast=\max\limits_{1\le j\le m}\{l_{j,j}\}, \quad q_j=\max\limits_{j\le i\le m} \{l_{i,j}\}.
   \end{equation}

Recall a homogeneous differential operator $A(D)=\sum\limits_{|\alpha|= \ell}a_{\alpha} D^\alpha$ is an elliptic if $A(\xi)=\sum\limits_{|\alpha|= \ell}a_\alpha \xi^\alpha>0,\, \,  \forall\xi\in\R^N,\, \,  \xi\neq 0$.

Further, if $B$ is a Banach space, then we denote by ${\mathbf B}$ the $m$-times topological direct product  
\[
{\mathbf B}  = B \otimes \dots \otimes B,
\] 
of spaces $B$. Elements of ${\mathbf B}$ are vector-functions $\Phi (x) = \langle{ \varphi_1(x), \dots, \varphi_m(x)} \rangle,$ where $\varphi_j(x) \in B, j=1,\dots,m.$ We also introduce the norm 
\[
||\Phi||^2_{\mathbf B}=||\varphi_1||_B^2+\cdots+||\varphi_m||_B^2\,.
\]

As usual, the symbol ${\mathbf C} ([0,T]; {\mathbf B})$ stands for the space of vector-functions continuous on $[0,T]$ with values in ${\mathbf B}$.

   The object of study in this work is the following problem
   
\textbf{Initial-Boundary Value Problem}. \textit{Find  functions $u_j(t,x)\in L_2^{q_j}(\mathbb{R}^n)$,
$t\in (0, T]$, $j=1, \dots, m$ (note that this inclusion is considered
as a boundary condition at infinity), such that}
\begin{equation}\label{conditionU}
 {U}(t,x)\in {\mathbf C}([0,T]\times \mathbb{R}^n
),  \, \,   
 D_t^{\mathcal{B}} {U}(t,x)\,\,\,\text{and} \,\,\, \mathbb{A}(D) {U} (t,x) \in {\mathbf C}((0,T]\times \mathbb{R}^n),   
\end{equation}
\textit{and satisfying the Cauchy problem}
\begin{align} \label{Cauchy_01_h}
D_t^{\mathcal{B}} {U}(t,x) &+ \mathbb{A}(D) {U} (t,x)=H(x,t), \quad t>0, \ x \in \mathbb{R}^n,
\\
U(0,x) & =\varPhi(x), \quad x \in \mathbb{R}^n, \label{Cauchy_02_h}
\end{align}
\textit{where $H(t,x)$ and $\varPhi(x)$ are given continuous functions.}

The solution to an initial-boundary value problem from class (\ref{conditionU}) is usually called \textit{a classical solution}.

Note that the requirement $u_j(t, x)\in L_2^{q_{j}}(\mathbb{R}^n)$ is necessary to prove the uniqueness of the solution, and on the other hand, the solution we found by the Fourier method, satisfies the specified condition.

We emphasize that the system of differential equations (\ref{Cauchy_01_h})  is parabolic in the sense of I. Petrovsky.
Let us recall the corresponding definition: a system (\ref{Cauchy_01_h}) is called parabolic in the sense of I. Petrovsky \cite{Petrowsky} if there is a positive constant $\delta$ such that $\Re(\mathcal{A}(\xi ) \mu, \mu)\geq \delta |\mu|^2$, for all $\mu\in \mathbb{C}^m$ (see also \cite{Kochuei1}).

The present paper consists of five sections. Section 2 provides some preliminary information about the Fourier and Laplace transforms, and the Mittag-Leffler functions. In Section 3 we construct a formal solution and in Section 4 we present results on the existence and uniqueness of a solution to Initial-Boundary Value Problem. The article ends with a “Conclusions” section, which discusses possible generalizations of the results obtained.
\section{Preliminaries}

\begin{definition}\label{Fourier_transform}

The Fourier transform of a function \( f: \mathbb{R}^n \rightarrow \mathbb{R} \), denoted as \( F[f](\xi) \) (or \( \hat{f}(\xi) \)), is defined as follows:
    \[
F[f](\xi)=\hat{f}(\xi) = \int\limits_{\mathbb{R}^n} f(x) e^{i x \xi} dx, \quad \xi \in \mathbb{R}^n.
\]
The inverse Fourier transform is denote by $F^{-1}[f](x)$ and is defined by the following integral:
\[
F^{-1}[f](x) = \frac{1}{(2 \pi)^n} \int\limits_{\mathbb{R}^n} f(\xi) e^{-i x \xi} d \xi, \quad x \in \mathbb{R}^n.
\]
\end{definition}

If, for example, $f\in L_2(\mathbb{R}^n)$, then $F[F^{-1}[f]](x)= F^{-1}[F[f]](x)=f(x)$.

If $f(x)$ is a function that is smooth enough and decays rapidly as $|x|$ goes to infinity, then the operation of the operator $A(D)$ on this function can be expressed as:
\begin{equation*}
A(D)f (x)= \frac{1}{(2 \pi)^n} \int\limits_{\mathbb{R}^n} A(\xi) F [f] (\xi)
e^{-i x \xi} d \xi\, \quad x \in \mathbb{R}^n.
\end{equation*}
Therefore for such a function one has
\begin{equation}\label{A(D)}
F[A(D) f(x)] (\xi)=A(\xi) F[f] (\xi).
\end{equation}

The classical Laplace transform is defined by the following integral formula:
$$
L[f](s)=\int\limits_{0}^{\infty}f(t)e^{-st}dt,\quad s\in \mathbb{C}.
$$
provided that the function $f$ (the Laplace original) is absolutely integrable on the semi-axis $\mathbb{R}_+$. In this case, the Laplace transform of a function, represented as \( L[f](s) \), is defined and analytic in the half-plane \(\Re s > 0\). Sometimes, the Laplace transform can be extended analytically to the left of the imaginary axis \(\Re s = 0\) into a larger region. This means there exists a non-positive real number \(\sigma_s\) such that \( L[f](s) \) remains analytic in the half-plane \(\Re s \ge \sigma_s \). Consequently, the following inverse Laplace transform can be introduced:
\[
L^{-1}[g](y)=\frac{1}{2\pi i}\int\limits_{\mathcal{L}_{ic}} e^{sy} L[f](s) ds.
\]
where $\mathcal{L}_{ic}=(c-i\infty,c+i\infty), c >\sigma_s$.

The relationship between the Laplace transform and its inverse is similar to that of the Fourier transform. If we apply the Laplace transform followed by the inverse Laplace transform to a "good enough" function $f(t)$, we should ideally obtain the original function back: $L[L^{-1}[f]](t)= L^{-1}[L[f]](t)=f(t)$.

The Laplace transform of the
Caputo derivative of a  function $f\in C^1(\mathbb{R}_+)$, 
is (see \cite{KLB}, p. 98):
\[
L[D_t^{\beta} f](s)= s^{\beta}L[{f}](s) -
f(0) s^{\beta-1},\,\,\, s>0.
\]
We will use this formula in the vector form.   For a vector function $\langle f_1, \dots, f_m \rangle$ we define:
\[
L[{D}_t^{\mathcal{B}} \langle f_1,\dots,f_m\rangle](s)=\langle L[{D}_t^{\beta_1}f_1](s), \dots,  L[{D}_t^{\beta_m} f_m](s)\rangle.
\]
Then
\begin{equation}\label{CD_LT_01} 
L[{D}_t^{\mathcal{B}} \langle f_1,\dots,f_m\rangle](s) = \langle s^{\beta_1} L[f_1](s) - f_1(0) s^{\beta_1-1},\dots,  s^{\beta_m} L[f_m](s)  - f_m(0) s^{\beta_m-1}\rangle.
\end{equation}

Additionally we use the following formula (see \cite{KLB}, p. 84):
$$
L[{\partial}_t^{\mathcal{B}} \langle f_1,\dots,f_m\rangle](s) = \langle L[{\partial}_t^{\beta_1}f_1](s), \dots,  L[{\partial}_t^{\beta_m} f_m](s)\rangle=
$$
\begin{equation}\label{R-LLT}
    =\langle s^{\beta_1} L[f_1](s) - I_t^{1-\beta_1}(f_1)(0+),\dots,  s^{\beta_m} L[f_m](s)  - I_t^{1-\beta_m}(f_m)(0+) \rangle,
\end{equation}
where $\partial_t^{\beta_j}$ and $I_t^{\beta_j}$ are Riemann-Liouville fractional derivative and integral, respectively and defined as:
$$
I^{\beta_j}_t f(t)= \frac{1}{\Gamma(1-\beta_j)}\int\limits_0^t \frac{f(\xi)}{(t-\xi)^{1-\beta_j}}d\xi, \quad \partial^{\beta_j}_t f(t)= \frac{d}{dt}I^{1-\beta_j}_t f(t) , \,\,\, t>0.
$$

Recall that the Mittag-Leffler function $E_{\rho,\mu}(t)$ is defined as follows:
$$
E_{\beta,\mu}(t)=\sum_{k=0}^\infty \frac{t^k}{\Gamma(\rho k+\mu)},\quad \beta>0,\, \,  \mu\in\C.
$$
If the parameter $\mu =1$, then we have the classical
Mittag-Leffler function: $ E_{\beta}(z)= E_{\beta, 1}(z)$. Note also $E_{1, 1}(z)=E_{1}(z)=e^z$.

\begin{lemma}\label{MLest1}
    For the Mittag-Leffler function $E_\beta(-t)$ the following estimate is true (see, for example, \cite{Gor}, p. 174):
$$
 |E_\beta(-t)|\le\frac{C}{1+t},\,\, t\geq 0.
$$
\end{lemma}
\begin{lemma}\label{MLest}
(see, for example \cite{AF})    For $\lambda>0$ and $\varepsilon\in(0,1)$, there exists a constant $C > 0$, such that:
    $$
    |t^{\beta-1}E_{\beta,\beta}(-\lambda t^\beta)|\le C\lambda^{\varepsilon-1}t^{\varepsilon\beta-1}, \,\, t>0.
    $$
\end{lemma}

\begin{lemma}\label{MLLT}(see, for example, \cite{KLB}, p. 50)
    The Laplace transform of the Mittag-Leffler function is as follows:
    $$
        L[t^{\beta-1}E_{\beta,\beta}(-\lambda t^\beta)](s)=\frac{s^{\rho-\mu}}{s^{\rho}+\lambda}.
    $$
\end{lemma}

The convolution of $f(t)$ and $g(t)$ is written as "$f\ast g$" and is defined as follows
$$
(f\ast g)(t)=\int\limits_0^t f(\tau)g(t-\tau)d\tau=\int\limits_0^t f(t-\tau)g(\tau)d\tau.
$$

Laplace transform of convolution equal to this:
\begin{equation}\label{astlap}
L[f\ast g](s)=L[f](s)L[g](s).
\end{equation}

\section{Formal solution of Initial-Boundary Value Problem}

In this section, following the work of S. Umarov \cite{Umarovappear}, we will construct a formal solution to the problem (\ref{Cauchy_01_h})-(\ref{Cauchy_02_h}). To make the reasoning more understandable, let us first consider the case $m=3$ and, having explained the main ideas, then move on to the general case $\forall m\in \mathbb{N}$.

\subsection{Case $m=3$}

We solve the problem (\ref{Cauchy_01_h})-(\ref{Cauchy_02_h}) when $H(t,x)\equiv 0$. Applying the Fourier transform we get the following problem:
$$
\left\{
\begin{aligned}
    & D_t^{\mathcal{B}} F[U](t,\xi) +\mathcal{A}(\xi) F[U](t,\xi)=0, \quad t>0, \ \xi \in \mathbb{R}^n,\\
&F[U](0,\xi) = F[\varPhi](\xi), \quad \xi \in \mathbb{R}^n.
\end{aligned}
\right.
$$

Then apply the Laplace transform to get
$$
\begin{bmatrix}
s^{\beta_1} + A_{1,1}(\xi)& 0 & 0  \\
A_{2,1}(\xi) & s^{\beta_2}+A_{2,2}(\xi)& 0\\
A_{3,1}(\xi)  & A_{3,2}(\xi) & s^{\beta_3}+A_{3,3} (\xi)
\end{bmatrix}
\begin{bmatrix}
	L[F[u_1]](s,\xi)\\
	L[F[u_2]](s,\xi)\\
	L[F[u_3]](s,\xi)\\
\end{bmatrix}
$$
$$
=
\begin{bmatrix}
	s^{\beta_1-1}F[\varphi_1](\xi)\\
	s^{\beta_2-1}F[\varphi_2](\xi)\\
	s^{\beta_3-1}F[\varphi_3](\xi)\\
\end{bmatrix}.
$$
Solving this equations we get the following recurrent equality:
$$
	L[F[u_1]](s,\xi)=\frac{s^{\beta_1-1}}{s^{\beta_1} + A_{1,1}(\xi)}F[\varphi_1](\xi),
$$
$$
L[F[u_k]](s,\xi)=\frac{s^{\beta_k-1}}{s^{\beta_k} + A_{k,k}(\xi)}F[\varphi_k](\xi)	-\sum\limits_{j=1}^{k-1}\frac{A_{k,j}(\xi)}{s^{\beta_k} + A_{k,k}(\xi)}L[F[u_j]](s,\xi),\quad k=2,3
$$
From this we have:
$$
L[F[u_1]](s,\xi)=\frac{s^{\beta_1-1}}{s^{\beta_1} + A_{1,1}(\xi)}F[\varphi_1](\xi),
$$
$$
L[F[u_2]](s,\xi)=\frac{s^{\beta_2-1}}{s^{\beta_2} + A_{2,2}(\xi)}F[\varphi_2](\xi)-\frac{A_{2,1}(\xi)}{s^{\beta_2} + A_{2,2}(\xi)}L[F[u_1]](s,\xi)
$$
$$
=\frac{s^{\beta_2-1}}{s^{\beta_2} + A_{2,2}(\xi)}F[\varphi_2](\xi)-\frac{A_{2,1}(\xi)}{s^{\beta_2} + A_{2,2}(\xi)}\frac{s^{\beta_1-1}}{s^{\beta_1} + A_{1,1}(\xi)}F[\varphi_1](\xi),
$$
and
$$
L[F[u_3]](s,\xi)=\frac{s^{\beta_3-1}}{s^{\beta_3} + A_{3,3}(\xi)}F[\varphi_3](\xi)-\frac{A_{3,1}(\xi)}{s^{\beta_3} + A_{3,3}(\xi)}L[F[u_1]](s,\xi)
$$
$$
-\frac{A_{3,2}(\xi)}{s^{\beta_3} + A_{3,3}(\xi)}L[F[u_2]](s,\xi)
=\frac{s^{\beta_3-1}}{s^{\beta_3} + A_{3,3}(\xi)}F[\varphi_3](\xi)-\frac{A_{3,1}(\xi)}{s^{\beta_3} + A_{3,3}(\xi)}\frac{s^{\beta_1-1}}{s^{\beta_1} + A_{1,1}(\xi)}F[\varphi_1](\xi)
$$
$$
+\frac{A_{2,1}(\xi)A_{3,2}(\xi)}{(s^{\beta_2} + A_{2,2}(\xi))(s^{\beta_3} + A_{3,3}(\xi))}\frac{s^{\beta_1-1}}{s^{\beta_1} + A_{1,1}(\xi)}F[\varphi_1](\xi)-\frac{A_{3,2}(\xi)}{s^{\beta_3} + A_{3,3}(\xi)}\frac{s^{\beta_2-1}}{s^{\beta_2} + A_{2,2}(\xi)}F[\varphi_2](\xi),
$$

Now we use the inverse Laplace transform. Using Lemma \ref{MLLT} and equality (\ref{astlap}) we get:
$$
F[u_1](t,\xi)=E_{\beta_1}(-A_{1,1}(\xi)t^{\beta_1})F[\varphi_1](\xi),
$$
$$
F[u_2](t,\xi)=E_{\beta_2}(-A_{2,2}(\xi)t^{\beta_2})F[\varphi_2](\xi)
$$
$$
-A_{2,1}(\xi)\left[E_{\beta_1}(-A_{1,1}(\xi)t^{\beta_1})\ast (t^{\beta_2-1}E_{\beta_2,\beta_2}(-A_{2,2}(\xi)t^{\beta_2}))\right]F[\varphi_1](\xi),
$$
and
$$
F[u_3](t,\xi)=E_{\beta_3}(-A_{3,3}(\xi)t^{\beta_3})F[\varphi_3](\xi)
$$
$$
-\left(A_{3,1}(\xi)\left[E_{\beta_1}(-A_{1,1}(\xi)t^{\beta_1})\ast (t^{\beta_2-1}E_{\beta_2,\beta_2}(-A_{2,2}(\xi)t^{\beta_2}))\right]\right.
$$
$$
-A_{2,1}(\xi)A_{3,2}(\xi)\left[E_{\beta_1}(-A_{1,1}(\xi)t^{\beta_1})\ast (t^{\beta_2-1}E_{\beta_2,\beta_2}(-A_{2,2}(\xi)t^{\beta_2}))\right. 
$$
$$
\left.\left.\ast (t^{\beta_3-1}E_{\beta_3,\beta_3}(-A_{3,3}(\xi)t^{\beta_3}))\right] \right)F[\varphi_1](\xi)
$$
$$
-A_{3,2}(\xi)\left[E_{\beta_2}(-A_{2,2}(\xi)t^{\beta_2})\ast (t^{\beta_3-1}E_{\beta_3,\beta_3}(-A_{3,3}(\xi)t^{\beta_3}))\right]F[\varphi_2](\xi).
$$
Using the inverse Fourier transform we have:
$$
u_1(t,x)=\frac{1}{(2\pi)^n}\int\limits_{\R^n}E_{\beta_1}(-A_{1,1}(\xi)t^{\beta_1})F[\varphi_1](\xi)e^{-ix\xi}d\xi, 
$$
$$
u_2(t,x)=\frac{1}{(2\pi)^n}\int\limits_{\R^n}E_{\beta_2}(-A_{2,2}(\xi)t^{\beta_2})F[\varphi_2](\xi)e^{-ix\xi}d\xi
$$
\begin{equation}\label{m=3u_2}
-\frac{1}{(2\pi)^n}\int\limits_{\R^n} A_{2,1}(\xi)\left[E_{\beta_1}(-A_{1,1}(\xi)t^{\beta_1})\ast (t^{\beta_2-1}E_{\beta_2,\beta_2}(-A_{2,2}(\xi)t^{\beta_2}))\right]F[\varphi_1](\xi)e^{-ix\xi}d\xi,
\end{equation}
and
$$
u_3(t,x)=\frac{1}{(2\pi)^n}\int\limits_{\R^n}E_{\beta_3}(-A_{3,3}(\xi)t^{\beta_3})F[\varphi_3](\xi)e^{-ix\xi}d\xi
$$
$$
-\frac{1}{(2\pi)^n}\int\limits_{\R^n} A_{3,2}(\xi)\left[E_{\beta_2}(-A_{2,2}(\xi)t^{\beta_2})\ast (t^{\beta_3-1}E_{\beta_3,\beta_3}(-A_{3,3}(\xi)t^{\beta_3}))\right]F[\varphi_1](\xi)e^{-ix\xi}d\xi
$$
$$
-\int\limits_{\R^n}\left(A_{3,1}(\xi)\left[E_{\beta_1}(-A_{1,1}(\xi)t^{\beta_1})\ast (t^{\beta_2-1}E_{\beta_2,\beta_2}(-A_{2,2}(\xi)t^{\beta_2}))\right]\right.
$$
$$
-\left.A_{2,1}(\xi)A_{3,2}(\xi)\left[E_{\beta_1}(-A_{1,1}(\xi)t^{\beta_1})\ast (t^{\beta_2-1}E_{\beta_2,\beta_2}(-A_{2,2}(\xi)t^{\beta_2})) \right.\right.
$$
\begin{equation}\label{m=3u_3}
\left.\left.\ast(t^{\beta_3-1}E_{\beta_3,\beta_3}(-A_{3,3}(\xi)t^{\beta_3}))\right] \right)F[\varphi_1](\xi)e^{-ix\xi}d\xi.
\end{equation}

Our main goal is to determine what order of polynomial in $\xi$ is multiplied by the expression $ F[\varphi_i](\xi),\,\, i = 1, 2, 3$. In particular,  in the case when $ m = 3 $, $ F[\varphi_1](\xi) $ is multiplied by polynomials of order $\ell_{2,1} + \ell_{3,2}$  and $\ell_{3,1}$, $ F[\varphi_2](\xi) $ is multiplied by a polynomial of order $\ell_{3,2}$, and $ F[\varphi_3](\xi) $ is multiplied by a polynomial of order $0$ in the expression $u_3(t,x)$.

 We now extend this analysis to the general case.

\subsection{Case $\forall m\in \mathbb{N}$}

Applying Fourier transform to (\ref{Cauchy_01_h}) and (\ref{Cauchy_02_h}) (when $H(t,x)\equiv 0$), we obtain a system of fractional order ordinary differential equation (see (\ref{A(D)})):
$$
    D_t^{\mathcal{B}} F[U](t,\xi) +\mathcal{A}(\xi) F[U](t,\xi)=0, \quad t>0, \ \xi \in \mathbb{R}^n,
$$
and the initial condition
$$
F[U](0,\xi) = F[\varPhi](\xi), \quad \xi \in \mathbb{R}^n.
$$

Apply the Laplace transform in the vector form \eqref{CD_LT_01}, to get
$$
\langle s^{\beta_1} L[F[u_1]](s,\xi), \dots, s^{\beta_m} L[F[u_m]](s,\xi)  \rangle + \mathcal{A}(\xi) L[F[U]](s,\xi) 
$$
$$
= \langle s^{\beta_1-1} F[\varphi_1](\xi), \dots, s^{\beta_m-1}F[\varphi_m](\xi)\rangle , \quad s>0, \ \xi \in \mathbb{R}^n.
$$
As in the case $m=3$, we write this equation in matrix form:
$$
\begin{bmatrix}
s^{\beta_1} + A_{1,1}(\xi)& 0 & \dots & 0  \\
A_{2,1}(\xi) & s^{\beta_2}+A_{2,2}(\xi)&\dots& 0\\
\dots   & \dots    & \dots & \dots  \\
A_{m,1}(\xi)  & A_{m,2}(\xi) &\dots & s^{\beta_m}+A_{m,m} (\xi)
\end{bmatrix}
\begin{bmatrix}
	L[F[u_1]](s,\xi)\\
	L[F[u_2]](s,\xi)\\
	\dots\\
	L[F[u_m]](s,\xi)\\
\end{bmatrix}
$$
$$
=
\begin{bmatrix}
	s^{\beta_1-1}F[\varphi_1](\xi)\\
	s^{\beta_2-1}F[\varphi_2](\xi)\\
	\dots\\
	s^{\beta_m-1}F[\varphi_m](\xi)\\
\end{bmatrix}.
$$

From this, we get the recurrent equation:

\begin{equation}\label{u_1yechimi}
	L[F[u_1]](s,\xi)=\frac{s^{\beta_1-1}}{s^{\beta_1} + A_{1,1}(\xi)}F[\varphi_1](\xi),
\end{equation}

\begin{equation}\label{yechimsystem1}
L[F[u_k]](s,\xi)=\frac{s^{\beta_k-1}}{s^{\beta_k} + A_{k,k}(\xi)}F[\varphi_k](\xi)	-\sum\limits_{j=1}^{k-1}\frac{A_{k,j}(\xi)}{s^{\beta_k} + A_{k,k}(\xi)}L[F[u_j]](s,\xi), \quad k=2,\dots, m.
\end{equation}

Before solving these equations, we consider the following sets. 

Let $ G_{k,j}=\{k, k-1, \dots, k-j\} $ be a set, where $1\le k\le m$ and $1\le j\le k-1 $. By $G_{k,j}^{(h)},\, \,   0\le h\le j-1$ we denote subsets of $G_{k,j}$ which formed by excluding $h$ elements except for $k$ and $k-j$ elements of the set. For example, let $k=6, j=3$. Then $G_{k,j}=\{6, 5, 4, 3\}$ and $G_{k,j}^{(0)}=G_{k,j}$, $G_{k, j}^{(1)}=\{\{6, 5, 3\}, \{6, 4, 3\}\}$ and $G_{k,j}^{(2)}=\{ 6, 3\}$.

From (\ref{u_1yechimi}) and the recurrence equation (\ref{yechimsystem1}), we have:
\begin{equation}\label{yechimsystem2}
	L[F[u_k]](s,\xi)=\frac{s^{\beta_k-1}}{s^{\beta_k} + A_{k,k}(\xi)}F[\varphi_k](\xi)	+\sum\limits_{j=1}^{k-1}\sum\limits_{r=0}^{j-1}\frac{P_{k,j,r}(\xi)s^{\beta_{k-j}-1}}{\prod\limits_{\tau\in G_{k,k-j}^{(r)}}( s^{\beta_\tau} + A_{\tau,\tau}(\xi))}F[\varphi_{k-j}](\xi).
\end{equation}
where $P_{k,j,r}(\xi), k=2, \dots, m, j=1, \dots, k-1,$ are combinations of multiplication and sum of functions $A_{\mu,\nu}(\xi), \mu=2, \dots, k, \nu=1, \dots, j-1$.

Applying the inverse Laplace transform to (\ref{u_1yechimi}) and (\ref{yechimsystem2}), we get:
$$
	F[u_1](t,\xi)=E_{\beta_1}(-A_{1,1}(\xi)t^{\beta_1})F[\varphi_1](\xi),
$$
$$
	F[u_k](t,\xi)=E_{\beta_k}(- A_{k,k}(\xi)t^{\beta_k})F[\varphi_k](\xi)
 $$
 $$
 +\sum\limits_{j=1}^{k-1}\sum\limits_{r=0}^{j-1}\left[P_{k,j,r}(\xi)\left(E_{\beta_{k-j}}(- A_{k-j,k-j}(\xi)t^{\beta_{k-j}})\ast V_{k,j,r}(t,\xi)\right)\right]F[\varphi_{k-j}](\xi),
$$
where (see formulas (\ref{m=3u_2}) and (\ref{m=3u_3}))
$$
V_{k,j,r}(t,\xi)=\ast \prod\limits_{\tau\in G_{k,k-j}^{(r)}\tau\neq k-j} t^{\beta_\tau-1}E_{\beta_\tau,\beta_\tau}(-A_{\tau,\tau}(\xi)t^{\beta_\tau}).
$$
Here "$\ast$" is the convolution operation, and "$\ast \prod$" is the
convolution product. Thus, the solution of the Cauchy problem (\ref{Cauchy_01_h})-(\ref{Cauchy_02_h}) has the representation when $H(t,x)\equiv 0$:
\begin{equation}\label{sysyechh=0}
U(t,x)=S(t,D)\Phi(x),\quad t>0,\, \,  x\in\R^n,
\end{equation}
where $S(t,D)$ stands for the pseudo-differential matrix operator with the matrix-symbol $S(t,\xi)$:
$$
    s_{k,j}(t,\xi)= \begin{cases} 
    0 & \text{if } k<j, \\ 
    E_{\beta_k}(-A_{k,k}(\xi)t^{\beta_k}) & \text{if } k=j, \\
    \sum\limits_{j=1}^{k-1}\sum\limits_{r=0}^{j-1}\left[P_{k,j,r}(\xi)\left(E_{\beta_{k-j}}(- A_{k-j,k-j}(\xi)t^{\beta_{k-j}})\ast V_{k,j,r}(t,\xi)\right)\right] & \text{if } k>j.
    \end{cases}
$$

From this we get the following formulas for solutions  $u_i(t,x), i=1,...,m$:
$$
u_1(t,x)=\frac{1}{(2\pi)^n}\int\limits_{\R^n} E_{\beta_1}(-A_{1,1}(\xi)t^{\beta_1})F[\varphi_1](\xi) e^{-ix\xi}d\xi,
$$
$$
u_k(t,x)=\frac{1}{(2\pi)^n}\int\limits_{\R^n} E_{\beta_k}(-A_{k,k}(\xi)t^{\beta_k})F[\varphi_k](\xi) e^{-ix\xi}d\xi
$$
$$
+\frac{1}{(2\pi)^n}\sum\limits_{j=1}^{k-1}\sum\limits_{r=0}^{j-1}\int\limits_{\R^n} \left[P_{k,j,r}(\xi)E_{\beta_{k-j}}(- A_{k-j,k-j}(\xi)t^{\beta_{k-j}})\ast V_{k,j,r}(t,\xi)\right]F[\varphi_{k-j}](\xi) e^{-ix\xi}d\xi.
$$

If $H(t,x) \not\equiv 0$, then the representation of a formal solution to problem (\ref{Cauchy_01_h})-(\ref{Cauchy_02_h}) is obtained from (\ref{sysyechh=0}) and the fractional Duhamel principle \cite{Umarov1, Umarov2}. The solution takes the following form:
$$
U(t,x)=S(t,D)\Phi(x)+\int\limits_0^t S(\eta,D)\partial_t^{1-\mathcal{B}}H(t-\eta,x)d\eta.
$$

Let us also introduce the pseudo-differential operator $S'(\eta,D)$, the entries of the matrix-symbol $S'(\eta,\xi)$, which have the form:
$$
    s'_{k,j}(\eta,\xi)= \begin{cases} 
    0 & \text{if } k<j, \\ 
    \eta^{\beta_k-1}E_{\beta_k, \beta_k}(-A_{k,k}(\xi)\eta^{\beta_k}) & \text{if } k=j, \\
    \sum\limits_{j=1}^{k-1}\sum\limits_{r=0}^{j-1}\left[P_{k,j,r}(\xi)\left(\eta^{\beta_{k-j}-1}E_{\beta_{k-j},\beta_{k-j}}(- A_{k-j,k-j}(\xi)\eta^{\beta_{k-j}})\ast V_{k,j,r}(\eta,\xi)\right)\right] & \text{if } k>j.
    \end{cases}
$$

\begin{lemma}\label{equalityrep}
    The following equality holds:
    $$
\int\limits_0^t S'(\eta,D)H(t-\eta,x)d\eta=\int\limits_0^t S(\eta,D)\partial_t^{{\bf 1}-\mathcal{B}}H(t-\eta,x)d\eta,
    $$
where ${\bf 1}$ is a $m$ dimentional vector: ${\bf 1}=\langle{1,..,1\rangle}$  and $\partial_t^{{\bf 1}-\mathcal{B}}$ is the Riemann-Liouville fractional derivative.    
\end{lemma}
\emph{Proof.} Using the Laplace transform we get (see (\ref{astlap}) and (\ref{R-LLT})):
$$
L[S'](s,D)L[H](s,D)
$$
$$
=\begin{cases} 
    0, & \text{if } k<j, \\ 
    \frac{1}{s^{\beta_k}+A_{k,k}(\xi)}L[H](s,D), & \text{if } k=j, \\
    \sum\limits_{j=1}^{k-1}\sum\limits_{r=0}^{j-1}\left[P_{k,j,r}(\xi)\left(\frac{1}{s^{\beta_1}+A_{1,1}(\xi)}\dots \frac{1}{s^{\beta_{k-1}}+A_{k-1,k-1}(\xi)}\right)\right]L[H](s,D), & \text{if } k>j.
    \end{cases}
$$
$$
L[S](s,D)L[\partial_t^{{\bf 1}-\mathcal{B}}H](s,D)
$$
$$
=\begin{cases} 
    0,  \\ 
    \frac{s^{\beta_k-1}}{s^{\beta_k}+A_{k,k}(\xi)}s^{1-\beta_k}L[H](s,D),  \\
    \sum\limits_{j=1}^{k-1}\sum\limits_{r=0}^{j-1}\left[P_{k,j,r}(\xi)\left(\frac{s^{\beta_{k-j}-1}}{s^{\beta_{k-j}}+A_{k-j,k-j}(\xi)}\frac{1}{s^{\beta_1}+A_{1,1}(\xi)}\dots \frac{1}{s^{\beta_{k-1}}+A_{k-1,k-1}(\xi)}\right)\right]s^{1-\beta_{k-j}}L[H](s,D),
    \end{cases}
$$
$$
=\begin{cases} 
    0, & \text{if } k<j, \\ 
    \frac{1}{s^{\beta_k}+A_{k,k}(\xi)}L[H](s,D), & \text{if } k=j, \\
    \sum\limits_{j=1}^{k-1}\sum\limits_{r=0}^{j-1}\left[P_{k,j,r}(\xi)\left(\frac{1}{s^{\beta_1}+A_{1,1}(\xi)}\dots \frac{1}{s^{\beta_{k-1}}+A_{k-1,k-1}(\xi)}\right)\right]L[H](s,D), & \text{if } k>j.
    \end{cases}
$$

From this we obtain:
$$
L[S'](s,D)L[H](s,D)=L[S](s,D)L[\partial_t^{{\bf 1}-\mathcal{B}}H](s,D).
$$

Hence
$$
\int\limits_0^t S'(\eta,D)H(t-\eta,x)d\eta=\int\limits_0^t S(\eta,D)\partial_t^{{\bf 1}-\mathcal{B}}H(t-\eta,x)d\eta.
$$
Lemma is proved.

Using Lemma \ref{equalityrep}, we can rewrite the formal solution of problem (\ref{Cauchy_01_h})-(\ref{Cauchy_02_h}) in the following form:
\begin{equation}\label{yechim}
        U(t,x)=S(t,D)\Phi(x)+\int\limits_0^t S'(\eta,D)H(t-\eta,x)d\eta,\quad 0\le t\le T, x\in\R^n.
    \end{equation}

\section{Main theorem}

Before presenting the main theorem, we show that the formal solution (\ref{yechim}) constructed in the previous section satisfies all the requirements of Initial-Boundary Value Problem. Let us start by checking the convergence of the expression $\mathbb{A}(D)U(t,x)$. Note that each $u_i(t,x), i=1,\dots, m$ is subject to operators from the column $i$ of the given matrix $\mathbb{A}(D)$. For example, the operators $A_{1,1}(D), A_{2,1}(D), \dots, A_{m,1}(D)$ act on $u_1(t,x)$.

Let us denote (see (\ref{yechim})):
$$
U(t,x)=W(t,x)+Y(t,x),
$$
where
$$
W(t,x)=S(t,D)\Phi(x),\quad Y(t,x)=\int\limits_0^t S'(\eta,D)H(t-\eta,x)d\eta, \quad 0 \le t \le T,\, \,   x \in \R^n.
$$
We estimate expressions $W(t,x)$ and $Y(t,x)$ separately. Let $w_i(t,x)$ and $y_i(t, x)$, where $i=1,\dots,m,$ be the entries of the vector functions $W(t, x)$ and $Y(t, x)$, respectively.

To investigate $w_i(t,x),\,\, i=1,\dots,m$, it is sufficient to consider the following integrals:
$$
w_{1,R}(t,x)=\int\limits_{|\xi|<R} E_{\beta_1}(-A_{1,1}(\xi)t^{\beta_1})\hat{\varphi}_1(\xi) e^{-ix\xi}d\xi,
$$
$$
w_{k,R}(t,x)=\int\limits_{|\xi|<R} E_{\beta_k}(-A_{k,k}(\xi)t^{\beta_k})\hat{\varphi}_k(\xi) e^{-ix\xi}d\xi
$$
\begin{equation}\label{w_mRtx}
+\sum\limits_{j=1}^{k-1}\sum\limits_{r=0}^{j-1}\int\limits_{|\xi|<R} Q_{k,j,r}(t,\xi)\hat{\varphi}_{k-j}(\xi) e^{-ix\xi}d\xi,
\end{equation}
where $k=2,.., m$ and
$$
Q_{k,j,r}(t,\xi)=P_{k,j,r}(\xi)\left(E_{\beta_{k-j}}(- A_{k-j,k-j}(\xi)t^{\beta_{k-j}})\ast V_{k,j,r}(t,\xi)\right).
$$

We have the following statement.

\begin{lemma}\label{estAU}
    1) To estimate all coefficients in front of each of the $\hat{\varphi}_i(\xi), \, \,  i=1, \dots, m$, it is sufficient to estimate the expression $A_{q,q}(\xi)Q_{m,m-i,0}(t,\xi), \, \,  q=1,\dots m$.
    
    2) This estimation is as follows:
    \begin{equation}\label{Mainest5}
        |A_{q,q}(\xi)Q_{m,m-i,0}(t,\xi)|\le C_{\beta_{i},\beta_{i+1},\dots,\beta_m,\varepsilon}|\xi|^{p^\ast- 
 l_{i,i}+(m-i)\varepsilon}t^{\varepsilon\sum_{j=i+1}^m\beta_j-\beta_i},
    \end{equation}
    where $p^{\ast}$ defined in (\ref{p}) and $\varepsilon$ is an arbitrary number from the interval $(0,1)$.
\end{lemma}
\emph{Proof.} Let us first consider the coefficients in front of $\hat{\varphi}_1(\xi)$. To do this we use the method of mathematical induction. If $m=1$ and $m=2$, the checking is easy. Let $m=3$. Then the expression $A_{q,q}(\xi)Q_{3,2,r}(t,\xi), \, \,  q=1,2,3,$ in front of $\hat{\varphi}_1(\xi)$ is equal to:
$$
\sum\limits_{r=0}^1 A_{q,q}(\xi)Q_{3,2,r}(t,\xi)= A_{q,q}(\xi)Q_{3,2,1}(t,\xi)+A_{q,q}(\xi)Q_{3,2,0}(t,\xi)
$$
$$
=-A_{q,q}(\xi)A_{3,1}(\xi)\left[E_{\beta_{1}}(- A_{1,1}(\xi)t^{\beta_{1}})\ast \left(t^{\beta_3-1}E_{\beta_3,\beta_3}(-A_{3,3}(\xi)t^{\beta_3})\right)\right]
$$
$$
+A_{q,q}(\xi)A_{2,1}(\xi)A_{3,2}(\xi)\left[E_{\beta_{1}}(- A_{1,1}(\xi)t^{\beta_{1}})\ast \left(t^{\beta_2-1}E_{\beta_2,\beta_2}(-A_{2,2}(\xi)t^{\beta_2})\right) \right.
$$
$$
\left.\ast \left(t^{\beta_3-1}E_{\beta_3,\beta_3}(-A_{3,3}(\xi)t^{\beta_3})\right)\right].
$$
From this, it is evident that to estimate $\sum\limits_{r=0}^1 A_{q,q}(\xi)Q_{3,2,r}(t,\xi)$ it is only necessary to estimate $A_{q,q}(\xi)Q_{3,2,0}(t,\xi)$. Using (\ref{maxorder}), the estimation of \(A_{q,q}(\xi)Q_{3,2,0}(t,\xi)\) proceeds as follows:
$$
|A_{q,q}(\xi)Q_{3,2,0}(t,\xi)|=\left|A_{q,q}(\xi)A_{2,1}(\xi)A_{3,2}(\xi)\left[E_{\beta_{1}}(- A_{1,1}(\xi)t^{\beta_{1}})\ast \left(t^{\beta_2-1}E_{\beta_2,\beta_2}(-A_{2,2}(\xi)t^{\beta_2})\right) \right.\right. 
$$
$$
\left. \left. \ast \left(t^{\beta_3-1}E_{\beta_3,\beta_3}(-A_{3,3}(\xi)t^{\beta_3})\right)\right]\right|
\le |\xi|^{p^\ast +l_{2,2}+l_{3,3}}
$$
$$
\times\left|\int\limits_0^t\left[\int\limits_0^\sigma E_{\beta_{1}}(- A_{1,1}(\xi)\eta^{\beta_{1}})(\sigma-\eta)^{\beta_2-1}E_{\beta_2,\beta_2}(-A_{2,2}(\xi)(\sigma-\eta)^{\beta_2})d\eta\right.\right]
$$
$$
\left.\times(t-\sigma)^{\beta_3-1}E_{\beta_3,\beta_3}(-A_{3,3}(\xi)(t-\sigma)^{\beta_3})d\sigma\right|.
$$
Using Lemma \ref{MLest} we get:
$$
 |A_{q,q}(\xi)Q_{3,2,0}(t,\xi)|\le |\xi|^{p^\ast+l_{2,2}+l_{3,3}}\int\limits_0^t\int\limits_0^\sigma |\xi|^{-l_{1,1}}\eta^{-\beta_1}(\sigma-\eta)^{\varepsilon\beta_2-1}|\xi|^{\varepsilon-l_{2,2}}d\eta(t-\sigma)^{\varepsilon\beta_3-1}|\xi|^{\varepsilon-l_{3,3}}d\sigma
 $$
 \begin{equation}\label{estQ320}
 \le C_{\beta_1,\beta_2}|\xi|^{p^\ast-l_{1,1}+2\varepsilon} \int\limits_0^t \sigma^{\varepsilon\beta_2-\beta_1}(t-\sigma)^{\varepsilon\beta_3-1}d\sigma \le C_{\beta_1,\beta_2,\beta_3,\varepsilon}|\xi|^{p^\ast-l_{1,1}+2\varepsilon}t^{\varepsilon(\beta_2+\beta_3)-\beta_1}.
\end{equation}

Now we demonstrate the validity of this estimation for $A_{q,q}(\xi)Q_{2,1,0}(t,\xi), \, \,  q=1,2$, in front of $\hat{\varphi}_1(\xi)$ in $w_{2,R}(t, x)$. The expression $A_{q,q}(\xi)Q_{2,1,0}(t,\xi)$ can be represented as:
$$
A_{q,q}(\xi)Q_{2,1,0}(t,\xi)=-A_{q,q}(\xi)A_{2,1}(\xi)\left[E_{\beta_{1}}(- A_{1,1}(\xi)t^{\beta_{1}})\ast \left(t^{\beta_2-1}E_{\beta_2,\beta_2}(-A_{2,2}(\xi)t^{\beta_2})\right)\right].
$$
Following the estimation process described above, we obtain
\begin{equation}\label{estQ210}
|A_{q,q}(\xi)Q_{2,1,0}(t,\xi)|\le C_{\beta_2,\varepsilon}|\xi|^{p^\ast-l_{1,1}+\varepsilon}t^{\varepsilon\beta_2-\beta_1}.
\end{equation}
Comparing the estimates (\ref{estQ320}) and (\ref{estQ210}), we can see that to estimate the expressions $A_{q,q}(\xi)$ and $A_{q,q}(\xi)Q_{2,1,0}(t,\xi)$ in front of $\hat{\varphi_1}(\xi)$ in $w_{1,R}(t,x)$ and $w_{2,R}(t,x)$ respectively, it is sufficient to estimate the expression $A_{q,q}Q_{3,2,0}(t,\xi)$ in front of $\hat{\varphi_1}(\xi)$ in $w_{3,R}(t,x)$. Hence statement 1) of the lemma is proven for $m=3$ and for coefficients of $\hat{\varphi}_1(\xi)$ . 

Since the estimate (\ref{estQ320}) implies (\ref{Mainest5}), then statement 2) of the lemma is also proven for this case.

Now let us assume that the following estimate
$$
    |A_{q,q}(\xi)Q_{m-1,m-2,0}(t,\xi)|=|Z_{m-1}(\xi)||K_{m-1}(t,\xi)|=|A_{q,q}(\xi)A_{2,1}(\xi)A_{3,2}(\xi)\dots A_{m-1,m-2}(\xi)|
$$
\begin{equation}\label{umest1}
 \times\left|E_{\beta_{1}}(- A_{1,1}(\xi)t^{\beta_{1}})\ast V_{m-1,m-2,0}(t,\xi)\right|
 \le C_{\beta_1,\beta_2,\dots,\beta_{m-1},\varepsilon}|\xi|^{p^\ast-l_{1,1}+(m-2)\varepsilon}t^{\varepsilon\sum_{j=2}^{m-1}\beta_j-\beta_1},
\end{equation}
where 
$$
V_{m-1,m-2,0}(t,\xi)=\ast \prod\limits_{\tau\in G_{m-1,1}^{(0)}\tau\neq 1} t^{\beta_\tau-1}E_{\beta_\tau,\beta_\tau}(-A_{\tau,\tau}(\xi)t^{\beta_\tau}),
$$
and
$$
Z_{m-1}(\xi)=A_{q,q}(\xi)A_{2,1}(\xi)A_{3,2}(\xi)\dots A_{m-1,m-2}(\xi),\quad K_{m-1}(t,\xi)=E_{\beta_{1}}(- A_{1,1}(\xi)t^{\beta_{1}})\ast V_{m-1,m-2,0}(t,\xi),
$$
is valid for $m-1$ and prove the corresponding estimate for $m$.

Let us first estimate the expression $A_{q,q}(\xi)Q_{m,m-1,0}(t,\xi)$ in front of $\hat{\varphi}_1(\xi)$ in the function $w_{m,R}(t,x)$ (see (\ref{w_mRtx})). We have
$$
|A_{q,q}(\xi)Q_{m,m-1,0}(t,\xi)|=|Z_{m-1}(\xi) A_{m,m-1}(\xi)|\times
 \left|K_{m-1}(t,\xi)\ast t^{\beta_m-1}E_{\beta_m,\beta_m}(-A_{m,m}(\xi)t^{\beta_m})\right|.
$$
Using (\ref{maxorder}) and applying estimate (\ref{umest1}), we obtain:
$$
|A_{q,q}(\xi)Q_{m,m-1,0}(t,\xi)|\le C_{\beta_1,\beta_2,\dots,\beta_{m-1},\varepsilon}|\xi|^{l_{m,m}+p^\ast-l_{1,1}+(m-2)\varepsilon}
$$
$$
\times\int\limits_0^t \eta^{\varepsilon\sum_{j=2}^m\beta_j-\beta_1}(t-\eta)^{\beta_m-1}E_{\beta_m,\beta_m}(-A_{m,m}(\xi)(t-\eta)^{\beta_m}).
$$
Now using Lemma \ref{MLest} we get:
\begin{equation}\label{umes2}
|A_{q,q}(\xi)Q_{m,m-1,0}(t,\xi)|\le C_{\beta_1,\beta_2,\dots,\beta_{m},\varepsilon}|\xi|^{p^\ast-l_{1,1}+(m-1)\varepsilon}t^{\varepsilon\sum_{j=2}^{m}\beta_j-\beta_1}.
\end{equation}
The right-hand side of the estimate (\ref{umes2}) also gives an upper bound for the expressions before $\hat{\varphi}_1(\xi)$ in the function $w_{m,R}(t,x)$.
Indeed, first note, that in the term $Q_{m,m-1,0}(t,\xi)$, there is a product of elements below the diagonal, namely, a product of elements $A_{i,i-1}(\xi),\, \, i=2,\dots,m$. All the remaining terms  $Q_{m,m-1,r}(t,\xi),\, \, r=1,\dots, m-2,$ do not contain elements belonging to the same row (see (\ref{m=3u_3}) for the case $m=3$). Hence according to the assumption of mathematical induction, statement 1) is proved for the coefficients of $\hat{\varphi}_1(\xi)$. Note that the estimate (\ref{umes2}) implies (\ref{Mainest5}), thus the statement 2) is also proved.

In the same way, it can be shown that for all $\hat{\varphi}_i(\xi), \, \, i=2,\dots,m-1$, it is sufficient to estimate $Q_{m,m-i,0}(t,\xi)$ in front of $\hat{\varphi}_i(\xi), \, \, i=2,\dots,m-1$, in the expression for function $w_{m,R}(t,x)$ (see (\ref{w_mRtx})). In particular, for $m=3$ the function $\hat{\varphi}_2(\xi)$ in the expression $w_{3,R}(t,x)$ is multiplied by polynomials in $\xi$ of order $\ell_{3,2}$, and $\hat{\varphi}_2(\xi)$ in the expression $w_{2,R}(t,x)$ is multiplied by polynomials in $\xi$ of order $0$ (see ( \ref{m=3u_2}) and (\ref{m=3u_3})). This shows that the function $\hat{\varphi}_2(\xi)$ is affected by combinations of operators below the diagonal, with the exception of the operators in column 1 of the matrix $\mathbb{A}( D)$. 

In the case of arbitrary $m$, each function $\hat{\varphi}_i(\xi)$ is acted upon by operators below the diagonal of the matrix $\mathbb{A}(D )$ starting from column $i$. By removing columns and rows from the matrix $\mathbb{A}(D)$, up to index $i$, we obtain a matrix $\mathbb{A}_i(D)$ of dimension $(m-i)\times(m-i)$, similar to $\mathbb{A}(D)$. For $\hat{\varphi}_i(\xi)$ we can now repeat the same reasoning in the matrix $\mathbb{A}_i(D)$, with the help of which we estimated the expression multiplied by $\hat{\varphi}_1(\xi)$ in the matrix $\mathbb{A}(D)$. From here it is not difficult to obtain both statement 1) and statement 2) of the lemma.

Lemma is proved.

For all $w_{k,R}(t,x)$ we get the following estimate:
$$
    ||A_{q,q}(D)w_{k,R}(t,x)||_{C(\R^n)} = \bigg|\bigg|A_{q,q}(D)\int\limits_{|\xi|<R} E_{\beta_k}(-A_{k,k}(\xi)t^{\beta_k})\hat{\varphi}_k(\xi) e^{-ix\xi}d\xi
    $$
    $$
    +A_{q,q}(D) \sum\limits_{j=1}^{k-1}\sum\limits_{r=0}^{j-1}\int\limits_{|\xi|<R} Q_{k,j,r}(t,\xi)\hat{\varphi}_{k-j}(\xi) e^{-ix\xi}d\xi\bigg|\bigg|
$$
  $$  
\le \bigg|\bigg|A_{q,q}(D)\int\limits_{|\xi|<R} E_{\beta_k}(-A_{k,k}(\xi)t^{\beta_k})\hat{\varphi}_k(\xi) e^{-ix\xi}d\xi\bigg|\bigg|
$$
$$
+\bigg|\bigg|A_{q,q}(D) \sum\limits_{j=1}^{k-1}\sum\limits_{r=0}^{j-1}\int\limits_{|\xi|<R} Q_{k,j,r}(t,\xi)\hat{\varphi}_{k-j}(\xi) e^{-ix\xi}d\xi\bigg|\bigg|
$$
$$
\le\int\limits_{|\xi|<R} \bigg|A_{q,q}(\xi)E_{\beta_k}(-A_{k,k}(\xi)t^{\beta_k})\hat{\varphi}_k(\xi)\bigg|d\xi
$$
$$
+\sum\limits_{j=1}^{k-1}\sum\limits_{r=0}^{j-1}\int\limits_{|\xi|<R} \bigg|A_{q,q}(\xi) Q_{k,j,r}(t,\xi)\hat{\varphi}_{k-j}(\xi)\bigg|d\xi.
$$
Using Lemma \ref{estAU} we have:
$$
||A_{q,q}(D)w_{k,R}(t,x)||_{C(\R^n)}\le \sum\limits_{i=1}^k C_{\beta_i,\dots,\beta_k,\varepsilon}t^{\varepsilon\sum_{j=i+1}^{k}\beta_j-\beta_i} \int\limits_{\R^n} |\xi|^{p^\ast -l_{i,i}+(m-i)\varepsilon}|\hat{\varphi}_{i}(\xi)|d\xi.
$$
Now set $\tau > n/2$ and choose $\varepsilon > 0$ such that $\tau-(m-i)\varepsilon > n/2$. Then using the Holder
inequality we get:
$$
||A_{q,q}(D)w_{k,R}(t,x)||_{C(\R^n)}\le \sum\limits_{i=1}^k C_{\beta_i,\dots,\beta_k,\varepsilon}t^{\varepsilon\sum_{j=i+1}^{k}\beta_j-\beta_i}  \int\limits_{\R^n} |\xi|^{\tau+p^\ast -l_{i,i}}|\hat{\varphi}_{i}(\xi)||\xi|^{-\tau+(m-i)\varepsilon}d\xi
$$
$$
\le \sum\limits_{i=1}^k C_{\beta_i,\dots,\beta_k,\varepsilon,\tau}t^{\varepsilon\sum_{j=i+1}^{k}\beta_j-\beta_i} \bigg|\bigg|(1+|\xi|^2)^{\frac{1}{2}(\tau+p^\ast-l_{i,i})}\hat{\varphi}_{i}(\xi)\bigg|\bigg|_{L_2(\R^n)}
$$

\begin{equation}\label{Mainres1}
    \le \sum\limits_{i=1}^kC_{\beta_i,\dots,\beta_k,\varepsilon}t^{\varepsilon\sum_{j=i+1}^{m}\beta_j-\beta_i}||\varphi_i||_{L_2^{\tau+p^\ast-l_{i,i}}(\R^n)}.
\end{equation}

Let us now estimate expressions $y_i(t,x),\, \,  i=1,\dots, m $. To do this it is sufficient to consider the following functions:
$$
y_{1,R}(t,x)=\int\limits_0^t\int\limits_{|\xi|<R} \eta^{\beta_1-1}E_{\beta_1,\beta_1}(-A_{1,1}(\xi)\eta^{\beta_1})\hat{h}_1(t-\eta,\xi) e^{-ix\xi}d\xi d\eta,
$$
$$
y_{k,R}(t,x)=\int\limits_0^t\int\limits_{|\xi|<R} \eta^{\beta_k-1}E_{\beta_k,\beta_k}(-A_{k,k}(\xi)\eta^{\beta_k})\hat{h}_k(t-\eta,\xi) e^{-ix\xi}d\xi d\eta
$$
$$
+\sum\limits_{j=1}^{k-1}\sum\limits_{r=0}^{j-1}\int\limits_0^t\int\limits_{|\xi|<R} S_{k,j,r}(\eta,\xi)\hat{h}_{k-j}(t-\eta,\xi) e^{-ix\xi}d\xi d\eta,
$$
where $k=2,\dots,m$ and
$$
S_{k,j,r}(\eta,\xi)=P_{k,j,r}(\xi)\left(\eta^{\beta_{k-j}-1}E_{\beta_{k-j},\beta_{k-j}}(- A_{k-j,k-j}(\xi)\eta^{\beta_{k-j}})\ast V_{k,j,r}(\eta,\xi)\right).
$$

\begin{lemma}\label{h_itxest}
    1) To estimate all coefficients in front of each  $\hat{h}_i(t-\eta,\xi), \, \,  i=1, \dots, m$, it is sufficient to estimate only the expression $A_{q,q}(\xi)S_{m,m-i,0}(t,\xi), \, \,  q=1,\dots m$ in front of each $\hat{h}_i(t,\xi)$.
    
    2) This estimation is as follows:
    $$
        |A_{q,q}(\xi)S_{m,m-i,0}(\eta,\xi)|\le C_{\beta_{i},\beta_{i+1},\dots,\beta_m,\varepsilon}|\xi|^{p^\ast- 
 l_{i,i}+(m-i)\varepsilon}\eta^{\varepsilon\sum_{j=i}^m\beta_j},
    $$
     where $p^{\ast}$ defined in (\ref{p})  and $\varepsilon$ is an arbitrary number from the interval $(0,1)$.
\end{lemma}
This lemma is proved analogously as Lemma \ref{estAU}.

For all $y_k(t,x)$, $k=1,.., m $, we get the following estimate:
$$
\bigg|\bigg|A_{q,q}(D)y_{k,R}(t,x)\bigg|\bigg|_{C(\R^n)}=\bigg|\bigg| \int\limits_0^t\int\limits_{|\xi|<R} A_{q,q}(\xi)\eta^{\beta_k-1}E_{\beta_k,\beta_k}(-A_{k,k}(\xi)\eta^{\beta_k})\hat{h}_k(t-\eta,\xi) e^{-ix\xi}d\xi d\eta
$$
$$
+\sum\limits_{j=1}^{k-1}\sum\limits_{r=0}^{j-1}\int\limits_0^t\int\limits_{|\xi|<R} A_{q,q}(\xi)S_{k,j,r}(\eta,\xi)\hat{h}_{k-j}(t-\eta,\xi) e^{-ix\xi}d\xi d\eta\bigg|\bigg|
$$
$$
\le \int\limits_0^t\int\limits_{|\xi|<R} \bigg|A_{q,q}(\xi)\eta^{\beta_k-1}E_{\beta_k,\beta_k}(-A_{k,k}(\xi)\eta^{\beta_k})\hat{h}_k(t-\eta,\xi) e^{-ix\xi}\bigg|d\xi d\eta
$$
$$
+\sum\limits_{j=1}^{k-1}\sum\limits_{r=0}^{j-1}\int\limits_0^t\int\limits_{|\xi|<R} \bigg|A_{q,q}(\xi)S_{k,j,r}(\eta,\xi)\hat{h}_{k-j}(t-\eta,\xi) e^{-ix\xi}\bigg|d\xi d\eta.
$$
Using Lemma \ref{h_itxest} we have:
$$
\bigg|\bigg|A_{q,q}(D)y_{k,R}(t,x)\bigg|\bigg|_{C(\R^n)}\le \sum\limits_{i=1}^k C_{\beta_i,\dots,\beta_k,\varepsilon}\int\limits_0^t\eta^{\varepsilon\sum_{j=i+1}^{k}\beta_j} \int\limits_{\R^n} |\xi|^{p^\ast -l_{i,i}+(m-i)\varepsilon}|\hat{h}_{i}(t-\eta,\xi)|d\xi d\eta.
$$
Now set $\tau > n/2$ and choose $\varepsilon > 0$ such that $\tau-(m-i)\varepsilon > n/2$. Then using the Holder
inequality we get:
$$
||A_{q,q}(D)y_{k,R}(t,x)||_{C(\R^n)}\le \sum\limits_{i=1}^k C_{\beta_i,\dots,\beta_k,\varepsilon}\int\limits_0^t\eta^{\varepsilon\sum_{j=i+1}^{k}\beta_j}  \int\limits_{\R^n} |\xi|^{\tau+p^\ast -l_{i,i}}|\hat{h}_{i}(t-\eta,\xi)||\xi|^{-\tau+(m-i)\varepsilon}d\xi d\eta
$$
$$
\le \sum\limits_{i=1}^k C_{\beta_i,\dots,\beta_k,\varepsilon,\tau} \max\limits_t\bigg|\bigg|(1+|\xi|^2)^{\frac{1}{2}(\tau+p^\ast-l_{i,i})}\hat{h}_{i}(t,\xi)\bigg|\bigg|_{L_2(\R^n)}
$$

\begin{equation}\label{Mainres2}
    \le \sum\limits_{i=1}^kC_{\beta_i,\dots,\beta_k,\varepsilon,\tau}\max\limits_{t\in[0,T]}||h_i(t,\cdot)||_{L_2^{\tau+p^\ast-l_{i,i}}(\R^n)}.
\end{equation}

Based on the proven estimates (\ref{Mainres1}), (\ref{Mainres2}), the following main result of the work is proved by verbatim repetition of the proof of Theorem 2 of the work \cite{AO}.

\begin{theorem}\label{Main_theorem}
    Let $\tau>\frac{n}{2}$ and $\varphi_i(x)\in L_2^{\tau+p^\ast-l_{i,i}}(\R^n)$, $i=1,\dots, m$ and $h_i(t,x)\in L_2^{\tau+p^\ast-l_{i,i}}(\R^n)$, $i=1,\dots, m$. Then the solution of Initial-Boundary Value Problem exists and is unique, and the solution of this problem has the representation (\ref{yechim}).
\end{theorem}

\begin{remark}\label{nessesary}
If all operators $A_{j,j}(D)$  have the same order, then $p^\ast-l_{i,i}=0$. Therefore, the sufficient condition for the initial functions and for the right-hand sides of the equations has the form $\varphi_j, \, h_j(t,x) \in L^{\tau}_2(\mathbb{R}^n),$ $\tau>\frac{n}{2}.$ Each such function, by virtue of the Sobolev embedding theorem, is continuous, which is required in the conditions of Initial-Boundary Value Problem. Note that, condition $\tau>\frac{n}{2}$ cannot be weakened, i.e. even if $\tau =\frac{n}{2}$ then in the class $L^{\tau}_2(\mathbb{R}^n)$ there are also unbounded functions.
\end{remark}

\begin{remark}
    If the given matrix $\mathbb{A}(D)$ is an upper triangular matrix, then a similar result to that in Theorem \ref{Main_theorem} holds; that is, under the same conditions on the initial function and on the right side of the equation, the statement of Theorem \ref{Main_theorem} remains valid. However, in this case, only the type of solution to the initial boundary value problem will change. Construction of a solution to the problem in this case should begin from the bottom, specifically by finding the function \(u_m(t,x)\).
\end{remark}

\section{Conclusions} 
The article is the first to consider a system of fractional partial differential equations in the case when the vector order of fractional derivatives $\mathcal{B}$ has different components $\beta_j\in (0,1]$ and $\beta_j$ are not necessarily rational numbers. It is assumed that the matrix of the system $\mathcal{A}(\xi)$ is lower triangular, i.e. elements above the diagonal are zero. Sufficient conditions on the initial function and the right-hand side of the equation are found to guarantee the existence and uniqueness of classical solutions to the Cauchy problem for such systems. In some cases these conditions are also necessary (see Remark \ref{nessesary}).

Of course, the most interesting continuation of this work is the removal of the condition that elements above the diagonal of the system matrix are equal to zero. To do this, apparently it is necessary to deal with case $m=2$. Next, having a representation of the solution, we can study a system of equations, the right-hand side of which depends nonlinearly on the desired solution to the system.

These generalizations are the subject of further research.

\section*{Acknowledgements}
The authors are grateful to Sh. A. Alimov for discussions of
these results. The authors acknowledge financial support from  the Ministry of Higher Education, Science and Innovative Development of the Republic of Uzbekistan, Grant No F-FA-2021-424.

\bibliographystyle{amsplain}

\end{document}